\documentclass[10pt]{article}

\usepackage[a4paper]{geometry}
\usepackage[LSF, T1]{fontenc}
\usepackage{indentfirst}
\usepackage{hyperref}
\usepackage{breakurl}
\usepackage{amssymb, amsmath, amsthm}
\usepackage{graphicx}
\usepackage[labelsep=none]{caption}
\usepackage[labelformat=simple]{subcaption}

\newcommand{\llbr}{[\![}
\newcommand{\rrbr}{]\!]}
\newcommand{\bfa}{\mathbf{a}}
\newcommand{\bfb}{\mathbf{b}}
\newcommand{\bfh}{\mathbf{h}}
\newcommand{\bfu}{\mathbf{u}}
\newcommand{\bfv}{\mathbf{v}}
\makeatletter \newcommand{\myknight}{\textrm{\usefont{LSF}{skaknew}{m}{n} \cfss@symknight}} \makeatother

\hypersetup{bookmarksnumbered, pdfstartview=FitH, pdfpagemode=UseNone}

\newcommand{\mynewtheorem}[2]{\newtheorem{#1}{\indent #2}}
\mynewtheorem{lemma}{Lemma}
\mynewtheorem{proposition}{Proposition}
\mynewtheorem{theorem}{Theorem}
\mynewtheorem{corollary}{Corollary}
\newenvironment{myproof}[1][Proof]{\begin{proof}[\indent #1]}{\end{proof}}

\begin{document}

\title{\textbf{Leaper Embeddings}}
\author{Nikolai Beluhov}
\date{}

\maketitle

\begin{abstract} A leaper is a chess piece which generalises the knight. Given $n$ and a $(p, q)$-leaper $L$, we study the greatest $m$ such that the $m \times m$ grid graph can be embedded into the $n \times n$ leaper graph of $L$. We can assume that $p$ and $q$ are relatively prime. We show that $m \approx n$ when $p$ and $q$ are of opposite parities and $m \approx n/2$ otherwise. The latter case is substantially more difficult. The proof involves certain combinatorial-geometric results on the chords of connected figures which might be of independent interest. \end{abstract}

\section{Introduction} \label{intro}

We say that a graph $G$ \emph{embeds} into a graph $H$ if $G$ is isomorphic to some subgraph $G'$ of $H$. We also call $G'$ an \emph{embedding} of $G$ into $H$.

Let $[n] = \{1, 2, \ldots, n\}$. We define the \emph{grid graph} $\square_n$ to be the graph on vertices $[n]^2$ where $(x', y')$ and $(x'', y'')$ are joined by an edge if they are at unit Euclidean distance apart. The \emph{knight graph} $\myknight_n$ uses the same vertex set but with the adjacency rule $\{|x' - x''|, |y' - y''|\} = \{1, 2\}$ instead. This is the (undirected) graph of all knight moves on an $n \times n$ chessboard.

Can we embed $\square_n$ into $\myknight_n$? If $n = 1$, then yes, trivially. Does any $n \neq 1$ admit an embedding? This is less obvious; in fact, the case of $n = 100$ was once posed as a mathematical olympiad problem. \cite{B} A little thought reveals that the answer is no -- there are no embeddings of $\square_n$ into $\myknight_n$ when $n \ge 2$.

\begin{figure}[ht] \centering \includegraphics{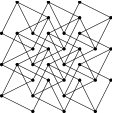} \caption{} \label{gk} \end{figure}

However, we can do almost as well as that. It is enough to make the grid just a bit smaller (or the knight graph just a bit bigger) for an embedding to become feasible. The solution to exercise 205 of \cite{K} shows that we can embed $\square_{n - 2}$ into $\myknight_n$ for all $n \ge 3$. For example, Figure \ref{gk} has $n = 8$.

One natural generalisation of knight graphs is given by leaper graphs. The $(p, q)$-\emph{leaper} $L$ is a fairy chess piece which moves similarly to the knight, but whose leap measures $p$ units along one coordinate axis and $q$ units along the other. We consider \emph{skew} leapers only, with nonzero and distinct $p$ and $q$. (The first few among them are the camel, $\{p, q\} = \{1, 3\}$; the giraffe, $\{p, q\} = \{1, 4\}$; and the zebra, $\{p, q\} = \{2, 3\}$.) We define the \emph{leaper graph} $\mathcal{L}_n$ to be the graph on vertices $[n]^2$ with adjacency rule $\{|x' - x''|, |y' - y''|\} = \{p, q\}$.

Our goal will be to estimate the greatest $m$ such that $\square_m$ embeds into a given leaper graph~$\mathcal{L}_n$.

First we note that a ``perfect'' embedding, with $m = n$, is possible only in the trivial case when $m = n = 1$. The argument of \cite{B} generalises immediately from the knight to all skew leapers -- see Section \ref{init}.

Observe next that, if $d = \gcd(p, q) \ge 2$, then the problem reduces from the $(p, q)$-leaper graph of order $n$ to the $(p/d, q/d)$-leaper graph of order $\lceil n/d \rceil$. So, from now on, we can safely assume that $\gcd(p, q) = 1$.

Since $\square_m$ is always connected, it must embed into a single connected component of $\mathcal{L}_n$. Let us, then, look into the way that $\mathcal{L}_n$ breaks down into connected components. When $\gcd(p, q) = 1$ and $p + q$ is odd, $L$ is known as \emph{free} and $\mathcal{L}_n$ is connected for all sufficiently large $n$. Otherwise, when $\gcd(p, q) = 1$ and both of $p$ and $q$ are odd, $L$ is known as \emph{half-free} and $\mathcal{L}_n$ consists of two connected components for all sufficiently large $n$. One of these connected components spans all vertices $(x, y)$ with $x + y$ even, and the other one spans those with $x + y$ odd.

We will show that all free leapers behave similarly to the knight; i.e., $m$ needs to be just slightly smaller than $n$ for an embedding to become feasible:

\begin{theorem} \label{tf} Suppose that $L$ is a free skew leaper. Then the greatest $m$ such that $\square_m$ embeds into $\mathcal{L}_n$ satisfies $m = n + O_L(1)$. \end{theorem} % thm free

(We use ``$O_L(1)$'' to mean ``a quantity whose absolute value is bounded from above by some constant which might depend on $L$ but certainly does not depend on $n$''.)

Since obviously $m \le n$, for the proof it suffices to find a construction which achieves the stated value of $m$. This is not too hard to do; we give one such construction in Section \ref{const}.

However, the difficulty jumps up sharply when we move on to half-free leapers. Then both connected components of $\mathcal{L}_n$ are of size about $n^2/2$, and so by way of an upper bound we get that $m$ cannot exceed approximately $n/\sqrt{2}$. On the other hand, the most natural ways to adapt the free-leaper construction to the half-free case only manage to achieve $m \approx n/2$; we give two such adaptations in Section \ref{const}. What is the true answer? We will show that the smaller value is in fact the right one:

\begin{theorem} \label{thf} Suppose that $L$ is a half-free skew leaper. Then the greatest $m$ such that $\square_m$ embeds into $\mathcal{L}_n$ satisfies $m = n/2 + O_L(1)$. \end{theorem} % thm half-free

Our proof involves some results in combinatorial geometry which might also be of independent interest outside of their connection to grid-and-leaper embeddings. We proceed now to a brief sketch of them.

We define a \emph{figure} in the plane to be any union of finitely many points, closed straight-line segments, and closed polygons. (So, in particular, our figures are always bounded.)

We say that a figure $F$ \emph{realises} a vector $\bfu$ if there exist two points of $F$ which differ by $\bfu$. Given two vectors $\bfu$ and $\bfv$, we also say that $\bfu$ \emph{forces} $\bfv$ if every connected figure which realises $\bfu$ must necessarily realise $\bfv$ as well.

When does $\bfu$ force $\bfv$? There is a classical result known as the ``universal chord theorem'' which states that, if $f$ is a continuous function over the real interval $[0; 1]$ with $f(0) = f(1)$, then for every positive integer $n$ we can find a horizontal chord in the graph of $f$ with length $1/n$. This generalises beyond the graphs of functions to arbitrary continuous curves -- see, for example, \cite{Y}. We get as a corollary of the generalisation that $\bfu$ forces $\bfv$ whenever $\bfu$ is an integer multiple of~$\bfv$.

We introduce the following variant of the forced chord problem: Let us say that the pair of vectors $\{\bfu_1, \bfu_2\}$ \emph{forces} the pair of vectors $\{\bfv_1, \bfv_2\}$ if every connected figure which realises both of $\bfu_1$ and $\bfu_2$ must necessarily realise at least one of $\bfv_1$ and $\bfv_2$ as well. When does one pair force another? We obtain one sufficient condition which is going to be crucial for our proof of Theorem~\ref{thf}. By way of an example, the simplest special case of this sufficient condition states that if $\bfu$ and $\bfv$ are non-parallel then $\{\bfu - \bfv, \bfu + \bfv\}$ forces $\{\bfu, \bfv\}$.

The rest of the paper is structured as follows: Section \ref{chord-i} is a quick review of the generalised universal chord theorem. In Section \ref{chord-ii}, we study the ``disjunctive'' variant of the forced chord problem, as outlined above. Section \ref{init} gathers some preliminary observations on grid-and-leaper embeddings. Then in Section \ref{const} we exhibit the aforementioned constructions which confirm both Theorem \ref{tf} and the lower bound of Theorem \ref{thf}. The upper bound of Theorem \ref{thf} is established in Section \ref{half}. Finally, in Section \ref{further} we collect some promising directions for further research.

\section{Forced Chords I} \label{chord-i}

Let $F$ be a connected figure and $\bfu$ a vector. Notice that $F$ realises $\bfu$ if and only if it realises~$-\bfu$. So changing the sign of a vector affects neither realisation nor forcing.

We write $F + \bfu$ for the translation copy of $F$ obtained when we translate it by $\bfu$. Clearly, $F$ realises $\bfu$ if and only if $F$ and $F + \bfu$ intersect.

Given two sets of vectors $U$ and $V$, we say that $U$ \emph{forces} $V$ if every connected figure which realises all elements of $U$ must necessarily realise at least one element of $V$ as well. We denote this relation between $U$ and $V$ by $U \vdash V$. When one of $U$ and $V$ is a singleton, we omit the curly braces since no confusion can arise.

Below (Lemma \ref{abc} and Proposition \ref{chord}) we reproduce the proof of the generalised universal chord theorem found in \cite{Y}, with some minor modifications.

\begin{lemma} \label{abc} Let $a$ and $b$ be positive real numbers with $a + b = 1$. Then $\bfu \vdash \{a\bfu, b\bfu\}$. \end{lemma}

\begin{myproof} Suppose, for the sake of contradiction, that $F$ realises $\bfu$ but not $a\bfu$ and $b\bfu$. Let $F_A = F - a\bfu$ and $F_B = F + b\bfu$. Let also $\ell'$ and $\ell''$ be the support lines of $F$ parallel to $\bfu$. Then $\ell'$ and $\ell''$ will be the support lines parallel to $\bfu$ of $F_A$ and $F_B$, too. (Figure \ref{abc-fig}.)

\begin{figure}[ht] \centering \includegraphics{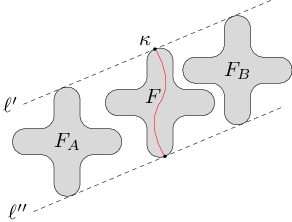} \caption{} \label{abc-fig} \end{figure}

Let $\kappa$ be any path in $F$ which connects two points on $\ell'$ and $\ell''$. Since $F$ does not realise $a\bfu$ and $b\bfu$, we get that $\kappa \subseteq F$ is disjoint from both of $F_A$ and $F_B$. Hence, $\kappa$ separates $F_A$ and $F_B$ within the strip bounded by $\ell'$ and $\ell''$. On the other hand, since $F$ realises $\bfu$, we get that $F_A$ and $F_B$ must intersect. This is a contradiction. \end{myproof}

\begin{proposition} \label{chord} Let $n$ be a positive integer. Then $\bfu \vdash \bfu/n$. \end{proposition}

\begin{myproof} By iterated application of Lemma \ref{abc}, if $a_1$, $a_2$, $\ldots$, $a_n$ are positive real numbers with $a_1 + a_2 + \cdots + a_n = 1$, then $\bfu \vdash \{a_1\bfu, a_2\bfu, \ldots, a_n\bfu\}$. The result follows when $a_1 = a_2 = \cdots = a_n = 1/n$. \end{myproof}

\section{Forced Chords II} \label{chord-ii}

We move on now to the ``disjunctive'' variant of the forced chord problem.

\begin{lemma} \label{abcd} Let $Q$ be a convex quadrilateral. Then the diagonals of $Q$ force the sides of $Q$. \end{lemma}

Or, more formally: Let $A$, $B$, $C$, $D$ be the vertices of $Q$, in this order. Then $\{A - C, B - D\} \vdash \{A - B, B - C, C - D, D - A\}$.

\begin{myproof} Let $F$ be a connected figure which realises both diagonals of $Q$. Let also $F_A$, $F_B$, $F_C$, $F_D$ be four translation copies of $F$ such that the pairwise displacements between them coincide with the pairwise displacements between $A$, $B$, $C$, $D$, respectively. (Figure \ref{abcd-fig}.)

\begin{figure}[ht] \centering \includegraphics{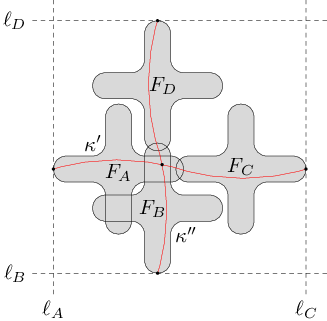} \caption{} \label{abcd-fig} \end{figure}

By means of a suitable affine transformation, we can assume without loss of generality that line $AC$ is horizontal with $A$ on the left and $C$ on the right; while line $BD$ is vertical with $B$ below and $D$ above. Let $\ell_A$, $\ell_B$, $\ell_C$, $\ell_D$ be, respectively, the left-hand vertical support line of $F_A$; the lower horizontal support line of $F_B$; the right-hand vertical support line of $F_C$; and the upper horizontal support line of $F_D$.

Then in fact all of $F_A$, $F_B$, $F_C$, $F_D$ will lie on the right of $\ell_A$; above $\ell_B$; on the left of $\ell_C$; and below $\ell_D$. Hence, $\ell_A$, $\ell_B$, $\ell_C$, $\ell_D$ enclose a rectangle $R$ which contains all four of $F_A$, $F_B$, $F_C$, $F_D$.

Since $F$ realises $A - C$, we get that $F_A$ and $F_C$ intersect. Let $\kappa'$ be a path in $F_A \cup F_C$ which connects two points on the vertical sides $\ell_A$ and $\ell_C$ of $R$. Similarly, let $\kappa''$ be a path in $F_B \cup F_D$ which connects two points on the horizontal sides $\ell_B$ and $\ell_D$ of $R$.

Since both of $\kappa'$ and $\kappa''$ are contained within $R$, they must intersect. Thus one of $F_A$ and $F_C$ must necessarily intersect one of $F_B$ and $F_D$. This means that $F$ does realise some side of $Q$, as needed. \end{myproof}

In Lemma \ref{abcd}, a pair forces a quadruple. By iterated application of this lemma, though, we can bring our set sizes down to a pair forcing a pair:

\begin{proposition} \label{fork} Let $a$, $b$, $c$, $d$ be four integers with $abcd < 0$. Let also $\bfu$ and $\bfv$ be two non-parallel vectors. Then $\{a\bfu + b\bfv, c\bfu + d\bfv\} \vdash \{\bfu, \bfv\}$. \end{proposition}

Setting $a = b = c = 1$ and $d = -1$, we obtain the sum-and-difference special case mentioned in the introduction. (The same special case also follows by Lemma \ref{abcd} directly, when we make $Q$ a parallelogram.)

\begin{myproof} By means of a suitable affine transformation, we can assume without loss of generality that $\bfu = (1, 0)$ and $\bfv = (0, 1)$. Let $F$ be any connected figure which realises both of $(a, b)$ and $(c, d)$. Our goal is to show that it realises at least one of $(1, 0)$ and $(0, 1)$, too.

We say that the vectors $(x, y)$ and $(-z, t)$ form a \emph{good} pair if $x$, $y$, $z$, $t$ are positive integers and $F$ realises both of them. For example, we can assume without loss of generality that the pair $\{(a, b), (c, d)\}$ is good. (By swapping the two vectors and changing their signs if necessary.)

We introduce the following procedure: Suppose that $\{(x, y), (-z, t)\}$ is a good pair. Let $w = \min\{x, z\}$. Then consider the four points \[\begin{aligned} A &= (-w/2, 0) & B &= (w/2, 0)\\ C &= (x - w/2, y) & D &= (w/2 - z, t). \end{aligned}\]

Clearly, $ABCD$ is a convex quadrilateral. (Figure \ref{step}.) Since $F$ realises both of its diagonals, by Lemma \ref{abcd} it must also realise at least one of its sides -- say, $\sigma$. If $\sigma$ is either horizontal or vertical, then by Proposition \ref{chord} we get that $F$ realises either $(1, 0)$ or $(0, 1)$, respectively, and we are done.

\begin{figure}[ht] \centering \includegraphics{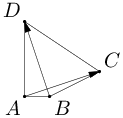} \caption{} \label{step} \end{figure}

Otherwise, the sign of the slope of $\sigma$ is well-defined. We discard whichever one of $(x, y)$ and $(-z, t)$ is of slope with the same sign, and we replace it with $\sigma$. Or, to be precise, with the appropriate difference between the endpoints of $\sigma$; we set up the sign so that the second coordinate of the replacement is positive.

This gives us a new good pair. Let us now iterate our procedure, beginning with $(a, b)$ and $(c, d)$. Since the process can only stop with an invocation of Proposition \ref{chord} which shows that $F$ realises one of $(1, 0)$ and $(0, 1)$, it suffices to demonstrate that we cannot keep going indefinitely.

Suppose that one iteration of our procedure has transformed the good pair $\{(x, y), (-z, t)\}$ into the good pair $\{(x', y'), (-z', t')\}$. (Of course, either $(x, y) = (x', y')$ or $(-z, t) = (-z', t')$.) We say that the iteration is \emph{regular} if one of $x$, $y$, $z$, $t$ has become smaller while the other three have remained the same. Or, in other words, if $x \ge x'$, $y \ge y'$, $z \ge z'$, $t \ge t'$ and exactly one of these inequalities is strict.

We claim that, over the course of the entire process, at most one irregular iteration can occur.

Indeed, if $\sigma = AB$, then $\sigma$ is horizontal and the process stops. If $\sigma = AD$, then either $\sigma$ is vertical and the process stops; or else $x < z$, we replace $(-z, t)$ with $(x - z, t)$, and the current iteration is regular. Similarly, if $\sigma = BC$, then either $\sigma$ is vertical and the process stops; or else $x > z$, we replace $(x, y)$ with $(x - z, y)$, and the current iteration is regular once again.

We are left to consider the case of $\sigma = CD$. If $y = t$, then $\sigma$ is horizontal and the process stops. Suppose, for concreteness, that $y < t$. (The analysis when $y > t$ is analogous.) So we must replace $(-z, t)$ with $(w - x - z, t - y)$. If $w = x$, the current iteration is regular. Otherwise, if $w = z < x$, the current iteration is irregular and we get a new good pair $\{(x', y'), (-z', t')\}$ with~$x' = z'$.

Or, to summarise, on an irregular iteration necessarily $\sigma = CD$ and in the new good pair we obtain as output the first coordinates of the two vectors share the same absolute value.

Suppose now that the current iteration has in fact been irregular and define $A'$, $B'$, $C'$, $D'$ as before but in terms of $\{(x', y'), (-z', t')\}$. Suppose also that the next iteration yields the good pair $\{(x'', y''), (-z'', t'')\}$. Since $A'B'$ is horizontal whereas both of $A'D'$ and $B'C'$ are vertical by virtue of $x' = z'$, we get that the realised side of $A'B'C'D'$ must have been $C'D'$. But then $x'' = z''$ once again. Furthermore, by looking into $y''$ and $t''$, too, we easily see that the next iteration must have been regular.

By induction, the same analysis will apply to all subsequent iterations as well. The property of the first coordinates of the two vectors in our good pair sharing the same absolute value will be maintained throughout; and all subsequent iterations will be regular.

This confirms our claim that at most one irregular iteration can take place over the course of the entire process. However, in the absence of irregular iterations, we can carry out only finitely many regular ones. The desired result now follows. \end{myproof}

We can also recast the continuous Propositions \ref{chord} and \ref{fork} in discrete form. Let $\square_\mathbb{Z}$ be the infinite grid graph on vertex set $\mathbb{Z}^2$. Consider any connected finite subgraph $G$ of $\square_\mathbb{Z}$. We say that $G$ realises $\bfu$ in the discrete sense if there exist two vertices of $G$ which differ by $\bfu$. The discrete forcing relation is defined as before, but based on the discrete notion of realisation.

Notice that we can view both grid and leaper graphs geometrically, with the vertices as integer points in the plane and the edges as straight-line segments joining these points. Let $G'$ be the figure associated with $G$. Clearly, $G$ realises a given integer vector $\bfu$ in the discrete sense if and only if $G'$ realises $\bfu$ in the continuous sense. Hence, if $U$ and $V$ are two sets of integer vectors such that $U$ forces $V$ in the continuous sense, then $U$ will also force $V$ in the discrete sense. The discrete forms of Propositions \ref{chord} and \ref{fork} are obtained when we apply this observation to the continuous originals. It is in fact these discrete forms that we are going to use in Section \ref{half}.

\section{Initial Observations on Embeddings} \label{init}

Let $\mathcal{L}_\mathbb{Z}$ be the infinite leaper graph of $L$ on vertex set $\mathbb{Z}^2$. Consider any embedding of $\square_m$ into $\mathcal{L}_\mathbb{Z}$, with $(i, j)$ mapping onto $P_{i, j}$ for all $i$ and $j$. Of course, every $4$-cycle of $\square_m$ must map onto a $4$-cycle of $\mathcal{L}_\mathbb{Z}$. But every $4$-cycle of $\mathcal{L}_\mathbb{Z}$ is (geometrically) a rhombus. So $P_{i, j} + P_{i + 1, j + 1} = P_{i, j + 1} + P_{i + 1, j}$ for all $i$ and $j$.

\begin{figure}[ht] \centering \includegraphics{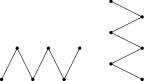} \caption{} \label{pk} \end{figure}

It follows that there exist $A_1$, $A_2$, $\ldots$, $A_m$ and $B_1$, $B_2$, $\ldots$, $B_m$ in $\mathbb{Z}^2$ with $P_{i, j} = A_i + B_j$ for all $i$ and $j$. Let $\alpha = A_1$---$A_2$---$\cdots$---$A_m$ and $\beta = B_1$---$B_2$---$\cdots$---$B_m$ be the two corresponding paths in $\mathcal{L}_\mathbb{Z}$. (Notice that, to us, a path is a certain kind of subgraph -- as opposed to a certain kind of a sequence of vertices within a graph.) Then we write $\alpha \times \beta$ for the original embedding of $\square_m$ into $\mathcal{L}_\mathbb{Z}$. For example, the $m = 6$ embedding in Figure \ref{gk} is the product of the two knight paths shown in Figure \ref{pk}.

It is straightforward now to rule out the possibility of a non-trivial perfect embedding. (This argument appears in \cite{B} for the special case of the knight.)

\begin{proposition} \label{imp} There is no embedding of $\square_n$ into $\mathcal{L}_n$ when $n \ge 2$. \end{proposition}

\begin{myproof} Suppose, for the sake of contradiction, that such an embedding does exist. Consider the four corner vertices of $\mathcal{L}_n$. Each one of them is of degree at most $2$. However, in $\square_n$ all non-corner vertices are of degree at least $3$. So the corner vertices of $\square_n$ must map onto the corner vertices of $\mathcal{L}_n$.

Observe next that the leap $A_1 \to A_2$ of $L$ points in the same direction as both of the leaps $P_{1, 1} \to P_{2, 1}$ and $P_{1, n} \to P_{2, n}$. Hence, from two distinct corner vertices of $\mathcal{L}_n$, there must exist two leaps of $L$ pointing in the same direction both of which land back within $\mathcal{L}_n$. This is a contradiction. \end{myproof}

Since the sums $A_i + B_j = P_{i, j}$ are pairwise distinct, we get that the pairwise differences between the vertices of $\alpha$ cannot overlap with the pairwise differences between the vertices of~$\beta$. Conversely, if $\alpha$ and $\beta$ are two paths in $\mathcal{L}_\mathbb{Z}$ that satisfy this condition, with $m$ vertices each, then there exists an embedding of $\square_m$ into $\mathcal{L}_\mathbb{Z}$ which factors as $\alpha \times \beta$.

However, we still want our embedding of $\square_m$ to fit within the subgraph $\mathcal{L}_n$ of $\mathcal{L}_\mathbb{Z}$. We go on to analyse this constraint as well.

Let $I$ and $J$ be two integer intervals (i.e., contiguous subsets of $\mathbb{Z}$) of sizes $k$ and $\ell$, respectively. We call $I \times J$ a \emph{box} of size $k \times \ell$. The \emph{bounding box} of a finite subset $S$ of $\mathbb{Z}^2$ is the smallest box which contains $S$; we denote it by $\llbr S \rrbr$. We also define the bounding box of a finite subgraph of $\square_\mathbb{Z}$ or $\mathcal{L}_\mathbb{Z}$ to be the bounding box of its vertex set. Clearly, if $\llbr \alpha \rrbr$ and $\llbr \beta \rrbr$ are of sizes $a_\mathsf{X} \times a_\mathsf{Y}$ and $b_\mathsf{X} \times b_\mathsf{Y}$, respectively, then $\llbr \alpha \times \beta \rrbr$ will be of size $(a_\mathsf{X} + b_\mathsf{X} - 1) \times (a_\mathsf{Y} + b_\mathsf{Y} - 1)$.

We are ready now to summarise our findings in the form of a handy criterion. (Similar observations appear in the solution to exercise 205 of \cite{K}, for the special case of the knight.)

\begin{lemma} \label{ab} There is an embedding of $\square_m$ into $\mathcal{L}_n$ if and only if there exist two paths $\alpha$ and $\beta$ in $\mathcal{L}_\mathbb{Z}$, with $m$ vertices each, such that: (i) The pairwise differences between the vertices of $\alpha$ are disjoint from the pairwise differences between the vertices of $\beta$; and (ii) The sizes $a_\mathsf{X} \times a_\mathsf{Y}$ of $\llbr \alpha \rrbr$ and $b_\mathsf{X} \times b_\mathsf{Y}$ of $\llbr \beta \rrbr$ satisfy $a_\mathsf{X} + b_\mathsf{X} \le n + 1$ and $a_\mathsf{Y} + b_\mathsf{Y} \le n + 1$. \end{lemma}

\section{Constructions} \label{const}

Suppose, to begin with, that $L$ is free. Below we present one construction which proves Theorem~\ref{tf}.

By swapping $p$ and $q$ if necessary, we can assume without loss of generality that $p$ is odd and $q$ is even. We write $\langle T \rangle^k$ for the concatenation of $k$ copies of the sequence $T$. Fix a positive integer $k$, and let $\alpha$ be the path traced out by $L$ when it makes $4kp$ successive leaps whose directions are specified by the sequence \[\big\langle \langle (p, q), (-p, q) \rangle^{p - 1}, (p, q), (p, q), \langle (p, -q), (-p, -q) \rangle^{p - 1}, (p, -q), (p, -q) \big\rangle^k.\]

\begin{figure}[ht] \centering \includegraphics{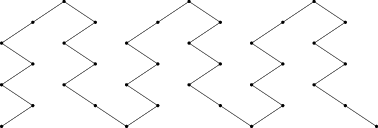} \caption{} \label{pz} \end{figure}

For example, Figure \ref{pz} shows $\alpha$ with $k = 3$ in the special case of the zebra, when $p = 3$ and~$q = 2$.

Let also $\beta$ be the $90^\circ$ counterclockwise rotation of $\alpha$. We proceed to check that $\alpha$ and $\beta$ satisfy the conditions of Lemma \ref{ab} with $m = 4kp + 1$ and $n = 4kp + 2pq + 1$. Of course, the desired result will then follow immediately.

\begin{figure}[ht] \centering \includegraphics{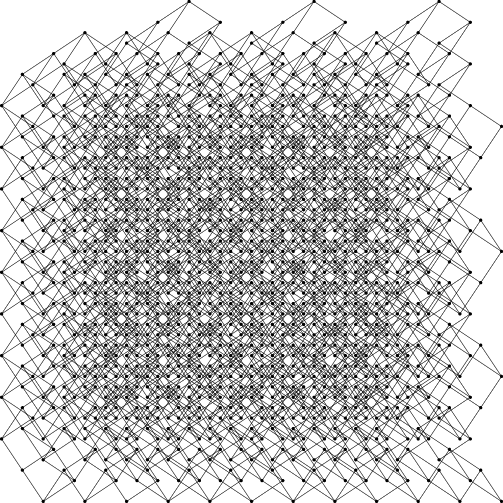} \caption{} \label{gz} \end{figure}

For example, Figure \ref{gz} shows $\alpha \times \beta$ with $\alpha$ as in Figure \ref{pz}; i.e., an embedding of the $37 \times 37$ grid graph into the $49 \times 49$ zebra graph.

We focus on part (i) of Lemma \ref{ab} as part (ii) is straightforward. Observe that each difference between two vertices of $\alpha$ is of the form $(px, qy)$ with $x + y$ even and $|y| \le 2p$. Similarly, each difference between two vertices of $\beta$ is of the form $(qz, pt)$ with $z + t$ even and $|z| \le 2p$.

Since $p \not \equiv q \pmod 2$ and $\gcd(p, q) = 1$, we get that the only potential overlaps are $(0, \pm 2pq)$, $(\pm 2pq, 0)$, and $(\pm 2pq, \pm 2pq)$. However, all differences $(px, qy)$ between the vertices of $\alpha$ with $y = \pm 2p$ satisfy $x \equiv 2 \pmod 4$. Since we broke symmetry in the beginning by ensuring that $p$ is odd and $q$ is even, it follows that $(0, \pm 2pq)$ and $(\pm 2pq, \pm 2pq)$ cannot occur as differences between the vertices of~$\alpha$. Similarly, $(\pm 2pq, 0)$ and $(\pm 2pq, \pm 2pq)$ cannot occur as differences between the vertices of~$\beta$. This completes the proof of Theorem \ref{tf}.

Suppose now that $L$ is half-free. We can derive the lower bound of Theorem \ref{thf} directly from the lower bound of Theorem \ref{tf}, as follows: Let $p' = |p - q|/2$ and $q' = (p + q)/2$. Then the $(p', q')$-leaper $L'$ is skew and free. Let also $d = \lfloor n/2 \rfloor$, and denote the $d \times d$ leaper graph of $L'$ by~$\mathcal{L}'_d$.

Consider the mapping $\varphi$ in the plane given by $(x, y) \to (x - y, x + y)$. (Geometrically, $\varphi$ is a spiral similarity centered at the origin, composed of a counterclockwise $45^\circ$ rotation and a dilation by a factor of $\sqrt{2}$.) Then $\varphi$ maps each leap of $L'$ onto a leap of $L$. Furthermore, the image $\varphi([d]^2)$ of the vertex set of $\mathcal{L}'_d$ is contained within a box of size $n \times n$. So $\mathcal{L}'_d$ embeds into~$\mathcal{L}_n$. Hence, the lower bound of Theorem \ref{tf} for $L'$ implies the lower bound of Theorem \ref{thf} for $L$.

Notice that any embedding of $\square_m$ into $\mathcal{L}_n$ obtained in this manner will be concentrated in the middle and surrounded on all sides by expanses of ``empty space''. Other approaches are possible as well where the embedding of $\square_m$ is spread out more evenly over the entirety of~$\mathcal{L}_n$. Below we present one such approach which proceeds by modifying the details of the free-leaper construction.

By swapping $p$ and $q$ if necessary, we can assume without loss of generality that $p \ge 3$. Fix a positive integer $k$, and let $\alpha$ be the path traced out by $L$ when it makes $2kp$ successive leaps whose directions are specified by the sequence \[\big\langle \langle (p, q), (-p, q) \rangle^{(p - 3)/2}, (p, q), (p, q), \langle (p, -q), (-p, -q) \rangle^{(p - 1)/2}, (p, -q), (p, q) \big\rangle^k.\]

Let also $\beta$ be the $90^\circ$ counterclockwise rotation of $\alpha$. We claim that $\alpha$ and $\beta$ satisfy the conditions of Lemma \ref{ab} with $m = 2kp + 1$ and $n = 4kp + pq + 1$. The proof follows along the same lines as the one in the free-leaper setting. The potential overlaps this time around work out to $(\pm pq, \pm pq)$. However, all differences $(px, qy)$ between the vertices of $\alpha$ with $y = p$ satisfy $x \equiv 3 \pmod 4$. Hence, $\alpha$ rules out half of the potential common differences. Similarly, $\beta$ rules out the other half. This confirms our claim. Just as before, the lower bound of Theorem~\ref{thf} is an immediate corollary.

\section{The Upper Bound in the Half-Free Case} \label{half}

We continue with the upper bound of Theorem \ref{thf}. Suppose, throughout this section, that $L$ is half-free. Let $\alpha$ and $\beta$ be two paths as in Lemma \ref{ab}, with $\llbr \alpha \rrbr$ of size $a_\mathsf{X} \times a_\mathsf{Y}$ and $\llbr \beta \rrbr$ of size $b_\mathsf{X} \times b_\mathsf{Y}$.

Each leap of $L$ is of slope either $\pm p/q$ or $\pm q/p$. By condition (i) of Lemma \ref{ab}, no such slope can occur in both of $\alpha$ and $\beta$. So we must split the four slopes between the two paths somehow.

The case when one of $\alpha$ and $\beta$ exhibits a single slope is straightforward. Indeed, suppose for concreteness that $p < q$ and all leaps of $L$ in $\alpha$ are of slope $p/q$. Then $m \le a_\mathsf{X}/q + O(1) \le n/3 + O(1)$ since $a_\mathsf{X} \le n$ by condition (ii) of Lemma \ref{ab} and $q \ge 3$.

Suppose, from now on, that $\alpha$ exhibits two out of the four possible slopes and $\beta$ exhibits the other two. This means that we can split the four vectors $(p, q)$, $(q, p)$, $(-p, q)$, $(-q, p)$ into two pairs $\{\bfa_\text{I}, \bfa_\text{II}\}$ and $\{\bfb_\text{I}, \bfb_\text{II}\}$ so that the directions of the leaps of $L$ in $\alpha$ are $\pm \bfa_\text{I}$ and $\pm \bfa_\text{II}$ whereas the directions of the leaps of $L$ in $\beta$ are $\pm \bfb_\text{I}$ and $\pm \bfb_\text{II}$.

Given two non-parallel vectors $\bfu$ and $\bfv$, we write $\mathfrak{L}(\bfu, \bfv)$ for the lattice generated by $\bfu$ and~$\bfv$. (Notice that, to us, a lattice is a certain kind of point set in the plane -- as opposed to a certain kind of subgroup within the additive group of all vectors in the plane.)

We also write $\mathcal{G}(\bfu, \bfv)$ for the graph on vertex set $\mathfrak{L}(\bfu, \bfv)$ where two vertices are joined by an edge if they differ by either $\bfu$ or $\bfv$. Clearly, this graph is always isomorphic to $\square_\mathbb{Z}$. In fact, geometrically, the two graphs are affine images of one another.

Let $\mathfrak{A} = \mathfrak{L}(\bfa_\text{I}, \bfa_\text{II})$ and $\mathfrak{B} = \mathfrak{L}(\bfb_\text{I}, \bfb_\text{II})$ as well as $\mathcal{A} = \mathcal{G}(\bfa_\text{I}, \bfa_\text{II})$ and $\mathcal{B} = \mathcal{G}(\bfb_\text{I}, \bfb_\text{II})$. By means of two suitable translations, we can assume without loss of generality that $\alpha$ is a path in $\mathcal{A}$ and $\beta$ is a path in $\mathcal{B}$.

Since both of $\mathcal{A}$ and $\mathcal{B}$ are affine images of $\square_\mathbb{Z}$, the discrete variants of Propositions \ref{chord} and~\ref{fork} apply to $\alpha$ and $\beta$. (Or, more precisely, to the appropriate affine images of $\alpha$ and $\beta$.) For convenience, from now on we will use the discrete notion of realisation when we talk about these paths; i.e., we are going to say that one of them realises some vector $\bfu$ if there exist two vertices of it which differ by $\bfu$. Condition (i) of Lemma \ref{ab} can now be restated as ``no nonzero vector is realised by both of $\alpha$ and $\beta$''.

Observe that the fundamental area of $\mathfrak{A}$ and $\mathfrak{B}$ is one and the same. (Explicitly, either $2pq$, $|p^2 - q^2|$, or $p^2 + q^2$, depending on how we split the slopes of $L$ between $\alpha$ and $\beta$.) Denote it by~$s$. Let also $h = s/2$ as well as $\bfh_\text{I} = (h, h)$ and $\bfh_\text{II} = (-h, h)$.

It is straightforward to see that $\mathfrak{A} \cap \mathfrak{B}$ is the lattice $\mathfrak{H} = \mathfrak{L}(\bfh_\text{I}, \bfh_\text{II})$. (This is true no matter how we split the slopes of $L$ between $\alpha$ and $\beta$, even though the details of the calculations are slightly different in each case.) So, in particular, $\pm \bfh_\text{I}$ and $\pm \bfh_\text{II}$ are the shortest nonzero vectors which could possibly be realised simultaneously by two paths in $\mathcal{A}$ and $\mathcal{B}$.

On the other hand, by condition (i) of Lemma \ref{ab}, neither one of $\bfh_\text{I}$ and $\bfh_\text{II}$ is realised by both of $\alpha$ and $\beta$. There are now two distinct cases to consider: Either $\alpha$ and $\beta$ miss one vector each, or else one of $\alpha$ and $\beta$ misses both vectors. These two cases are rather different, and we consider them separately below. The former case is easier, and in it the upper bound will follow by Proposition \ref{chord}. The latter case is tougher, and in it we are going to need Proposition \ref{fork} instead.

One simple observation will be useful in both cases: Since $\mathfrak{H}$ is a sub-lattice of $\mathfrak{A}$ and the fundamental area of $\mathfrak{H}$ is $2h^2 = hs$, we get that $\mathfrak{A}$ can be partitioned into $h$ translation copies of~$\mathfrak{H}$. Similarly, so can $\mathfrak{B}$. We denote $\mathcal{H} = \mathcal{G}(\bfh_\text{I}, \bfh_\text{II})$.

\smallskip

\emph{Case 1}. Each path misses one vector -- say, $\alpha$ misses $\bfh_\text{I}$ and $\beta$ misses $\bfh_\text{II}$.

Define a \emph{diagonal} to be any doubly infinite path in $\mathcal{H}$ all of whose edges are of the same slope; or any translation copy of such a path. Since $\alpha$ does not realise $\bfh_\text{I}$, by the discrete variant of Proposition \ref{chord} we get that $\alpha$ cannot contain two vertices in the same diagonal of slope $1$. Similarly, $\beta$ cannot contain two vertices in the same diagonal of slope $-1$.

\begin{figure}[ht] \centering \includegraphics{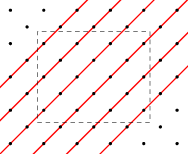} \caption{} \label{dia} \end{figure}

However, for each translation copy of $\mathfrak{H}$, its intersection with $\llbr \alpha \rrbr$ can be covered using $(a_\mathsf{X} + a_\mathsf{Y})/2h + O(1)$ diagonals of unit slope. (Figure \ref{dia}.) Since $\alpha$ meets each such diagonal at most once, and furthermore $\mathfrak{A}$ can be partitioned into $h$ translation copies of $\mathfrak{H}$, we find that $m \le (a_\mathsf{X} + a_\mathsf{Y})/2 + O_L(1)$. (The error term has increased from $O(1)$ to $O_L(1)$ because we are multiplying by $h$, which depends on $L$.) Similarly, also $m \le (b_\mathsf{X} + b_\mathsf{Y})/2 + O_L(1)$.

Finally, since $a_\mathsf{X} + a_\mathsf{Y} + b_\mathsf{X} + b_\mathsf{Y} \le 2n + 2$ by condition (ii) of Lemma \ref{ab}, we conclude that $m \le n/2 + O_L(1)$, as needed.

\smallskip

\emph{Case 2}. One path misses both vectors. Suppose, for concreteness, that $\alpha$ does.

Define a \emph{vertical zigzag} to be any doubly infinite path in $\mathcal{H}$ which can be oriented so that it alternates between moves in the directions $\bfh_\text{I}$ and $\bfh_\text{II}$; or any translation copy of such a path. Define also a \emph{horizontal zigzag} similarly, except that its moves must alternate between the directions $\bfh_\text{I}$ and $-\bfh_\text{II}$.

Suppose, for the sake of contradiction, that $\alpha$ contains two vertices in the same horizontal zigzag as well as two vertices in the same vertical zigzag. The former assumption tells us that $\alpha$ realises some vector of the form $c_1\bfh_\text{I} - c_2\bfh_\text{II}$, where $c_1$ and $c_2$ are nonnegative integers with $|c_1 - c_2| \le 1$. Similarly, the latter assumption tells us that $\alpha$ realises also some vector of the form $d_1\bfh_\text{I} + d_2\bfh_\text{II}$, where $d_1$ and $d_2$ are nonnegative integers with $|d_1 - d_2| \le 1$.

Since $\alpha$ does not realise $\bfh_\text{I}$ and $\bfh_\text{II}$, all four of $c_1$, $c_2$, $d_1$, $d_2$ are nonzero. However, this means that the discrete variant of Proposition \ref{fork} applies. By it, we obtain that $\alpha$ must realise at least one of $\bfh_\text{I}$ and $\bfh_\text{II}$ anyway. This is a contradiction.

Hence, either $\alpha$ meets each horizontal zigzag at most once or $\alpha$ meets each vertical zigzag at most once. Suppose, for concreteness, the latter.

\begin{figure}[ht] \centering \includegraphics{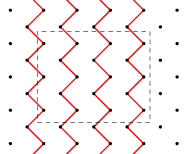} \caption{} \label{zz} \end{figure}

Observe now that, for each translation copy of $\mathfrak{H}$, its intersection with $\llbr \alpha \rrbr$ can be covered using $a_\mathsf{X}/2h + O(1)$ vertical zigzags. (Figure \ref{zz}.) Since $\alpha$ meets each such zigzag at most once, and furthermore $\mathfrak{A}$ can be partitioned into $h$ translation copies of $\mathfrak{H}$, we find that $m \le a_\mathsf{X}/2 + O_L(1)$. (Once again, the multiplication by $h$ has caused the error term to increase.)

Finally, since $a_\mathsf{X} \le n$ by condition (ii) of Lemma \ref{ab}, we conclude that $m \le n/2 + O_L(1)$, as needed. This completes the proof of Theorem \ref{thf}.

\section{Further Work} \label{further}

So far, in our study of forced chords we have focused on sufficient conditions. However, it would be interesting to consider necessity, too.

A well-known supplement to the universal chord theorem (see, for example, \cite{Y}) implies that, with $\bfu$ and $\bfv$ nonzero, $\bfu$ forces $\bfv$ if and only if $\bfu$ is an integer multiple of $\bfv$. Or, in other words, the sufficient condition of Proposition \ref{chord} is also necessary.

What about one set of vectors forcing another? Clearly, if $U$ forces $V$, also every superset of $U$ forces every superset of $V$. So it makes sense to restrict our analysis to the cases where neither $U$ nor $V$ can be replaced with a proper subset of itself. Then we say that $U$ \emph{minimally forces} $V$, and we write $U \models V$.

Suppose first that $U = \{\bfu\}$ is a singleton and $U \models V$. It is not too difficult to see that, in this setting, all elements of $V$ must be parallel to $\bfu$. Indeed, if a connected figure $F$ realises $\bfu$ but not any elements of $V$ parallel to $\bfu$, then by means of a suitable affine transformation we can ``squash'' $F$ into a new connected figure $F'$ which still realises $\bfu$ but misses $V$ altogether.

Hence, $V = \{a_1\bfu, a_2\bfu, \ldots, a_k\bfu\}$. Reasoning as in the proof of Proposition \ref{chord}, we see that $1$ being a linear combination of $|a_1|$, $|a_2|$, $\ldots$, $|a_k|$ with positive integer coefficients is a sufficient condition for $\{\bfu\}$ to force $V$. (Though perhaps not minimally.) Is it true that $\{\bfu\}$ minimally forces $V$ if and only if $a_1$, $a_2$, $\ldots$, $a_k$ satisfy this condition but no proper subset of them does?

Suppose now that $|U| = |V| = 2$ and the elements of $U$ and $V$ are pairwise non-parallel. Then is the sufficient condition for $U \vdash V$ of Proposition \ref{fork} also necessary?

Of course, these sub-problems are best viewed as stepping stones on the way to a general description of all instances in which one set of vectors forces another.

\section*{Acknowledgements}

The main results of this paper were obtained in 2021 after the author was contacted by Professor Donald Knuth in connection with exercises 205--209 of \cite{K}. The author is thankful to Prof.\ Knuth for bringing this series of exercises to his attention.

\end{document}